\documentclass[11pt]{article}
\usepackage{graphicx}
\usepackage{amssymb}
\usepackage{epstopdf}
\newcounter{example} %example counter

 \newcommand{\N}{{\bf N}}

 \newcommand{\Z}{{\bf Z}}

 \newcommand{\bigconv}{\ast \prod}
 \newcommand{\be}{\begin{equation}}
 \newcommand{\ee}{\end{equation}}
\begin{document}
\title{Cosine Products, Fourier Transforms, and Random Sums 
       \footnote{This article has appeared in {\em Amer. Math. Monthly}, 102:716--724, 1995.}}
\author{ Kent E. Morrison \\
\\Department of Mathematics 
\\California Polytechnic State University
\\San Luis Obispo, CA  93407
\\kmorriso@oboe.calpoly.edu}
\date{March 24, 1994}
\maketitle


%\begin{abstract}
%We investigate several infinite product of cosines and find the 
%closed form using the Fourier transform. The answers provide limiting 
%distributions for some elementary probability experiments. 
%\end{abstract} 
%\renewcommand{\baselinestretch}{2} \large \normalsize

\section{Introduction}
The function $\sin x /x$ is endlessly fascinating. By setting $x=\pi/2$ 
in the infinite product expansion
\begin{equation}
 \frac{\sin x}{x} = \prod_{k=1}^{\infty}\cos \frac{x}{2^k}
 \label{identity1}
\end{equation}
one gets the first actual formula for $\pi$ that mankind ever discovered,
dating from 1593 and due to
Fran\c{c}ois Vi\`ete (1540-1603), whose Latinized name is Vieta. 
(Was any notice taken of the formula's 400th anniversary, perhaps by the issue
of a postage stamp?) From the samples of a function $f(x)$ at equally spaced points $x_n$, $n \in 
\Z$, one can reconstruct the complete function with the aid of $\sin 
x/x$, provided $f$ is ``band-limited'' and the spacing of the samples is 
small enough. This is the content of the Sampling Theorem, which lends 
its name to $\sin x/x$ as the {\bf sampling function}. Its importance in 
signal processing, where it is also known as $\mbox{sinc}\ x$, is the result of its 
Fourier transform being the characteristic function of the interval
$[-1,1]$ (modulo a scalar factor).

In section 2 we prove the infinite product expansion for $\sin x/x$ and 
derive Vi\`ete's formula. In section 3 we transform the product expansion
with the Fourier transform and use convolution and delta distributions to 
prove it in a way that reveals a host of similar identities. Section 4 
puts these identities into a probabilistic setting, and in section 5 we 
alter the probability experiments in order to make connections between 
infinite cosine products, Cantor sets, and sums of series with random 
signs, particularly the harmonic series. This leaves us with some 
interesting unsolved problems and conjectures for further work.

\section{An Elementary Proof}
Repeated use of the double angle formula for the sine shows that
\begin{eqnarray*}
 \sin x & = & 2 \sin \frac{x}{2} \cos \frac{x}{2} \\
        & = & 4 \sin \frac{x}{4}  \cos \frac{x}{4} \cos \frac{x}{2} \\
        & \vdots &    \\
        & = & 2^n \sin \frac{x}{2^n}(\prod_{k=1}^n \cos \frac{x}{2^k}) 
\end{eqnarray*}
But 
\begin{displaymath}
	\lim_{n \rightarrow \infty} 2^n \sin \frac{x}{2^n} = x,
\end{displaymath} 
thereby proving the identity. See Figure 1 for an indication of how 
quickly the product converges.

Let $x=\pi/2$, make use of the half-angle identity, and there you have 
Vi\`ete's formula for $\pi$,
\begin{equation}
 \frac{2}{\pi}=\frac{\sqrt{2}}{2}\frac{\sqrt{2+\sqrt{2}}}{2}
               \frac{\sqrt{2+\sqrt{2+\sqrt{2}}}}{2}\cdots \, .
\end{equation}               
At this point the cosine identity could remain an isolated curiosity of 
historical interest, relegated to the ends of exercise sets in textbooks.
In fact, it is just the first of an infinite family of cosine product 
identities for $\sin x/x$. 

\vspace{.2in}
\centerline {
\includegraphics[width=3in]{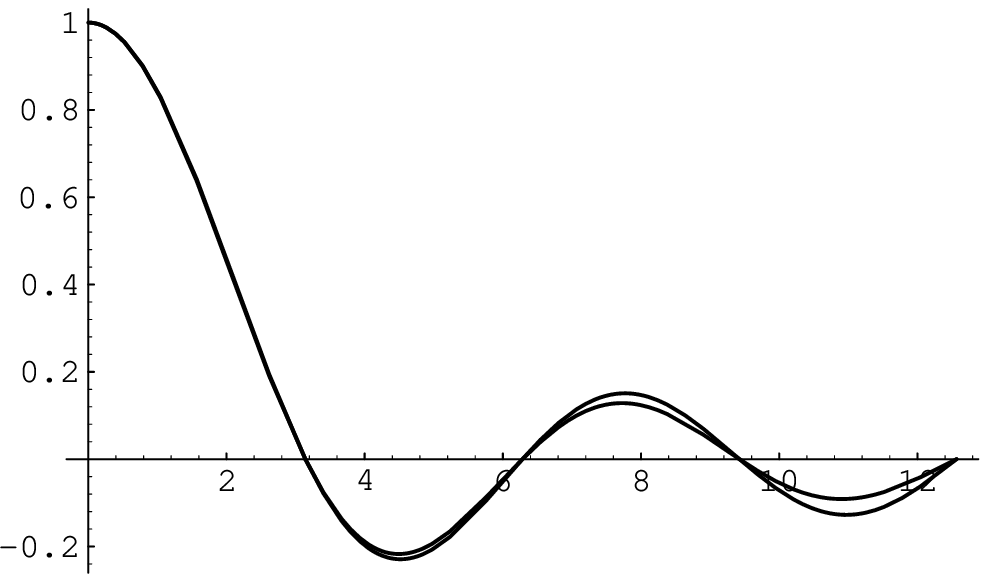}
}
\begin{center} 
	  Figure \thefigure. Graphs of $\sin x/x$ and 
	  $\cos \frac{x}{2} \cos \frac{x}{4} \cos \frac{x}{8}$. \\
	  Where both graphs
	  are visible, $\sin x/x$ is nearer the $x$ axis.  
\end{center}  
\addtocounter{figure}{1}

\section{The Fourier Transform and More Identities}

For a complex valued function $f(x)$ defined on the real line, the 
Fourier transform puts together $f$ as a continuous linear combination of 
the ``pure'' oscillations $e^{i\omega x}$ in which the coefficient in 
front of $e^{i\omega x}$ is denoted by $\hat{f}(\omega)$. Thus,
\begin{equation}
   f(x)= \int_{-\infty}^{\infty}\hat{f}(\omega) e^{i\omega x} \, d\omega .
\end{equation}
The function $\hat{f}$ is the {\bf Fourier transform} of $f$ and the 
integral above is a description of how to get back $f$ from $\hat{f}$
and is actually the formula for the inverse transform. How do we get 
$\hat{f}$ from $f$? That is given by this integral:
\begin{equation}
   \hat{f}(\omega)= \frac{1}{2\pi}\int_{-\infty}^{\infty}f(x)e^{-i\omega 
   x}\, dx .
\end{equation}
Of course, the proofs of these relationships involve hypotheses on the 
functions so that the integrals make sense, but they can be extended 
beyond the realm of ordinary functions to generalized functions or 
distributions. We need more than ordinary functions in order to make 
sense of the Fourier transform of a sine or cosine.

Notation: we also write the Fourier transform of $f$ as ${\cal F}(f)$ and 
the inverse transform of $\phi$ as ${\cal F}^{-1}(\phi)$. 

Consider $\cos bx$, which by Euler's Identity may be written as
\[ \cos bx = \frac{1}{2}(e^{ibx} + e^{-ibx} ).  \]
This shows the function written as a linear combination of just two of
of the functions $e^{i\omega x}$ for $\omega = b$ and $\omega = -b$. The 
coefficients appear to be $1/2$, but if we use them in the integral form 
with all other coefficients zero, then we cannot represent the cosine 
function. Instead, we must regard the coefficients as point masses at 
$b$ and $-b$. Therefore, the Fourier transform of $\cos bx$ is
$(1/2)(\delta_b + \delta_{-b})$, where $\delta_b$ denotes the Dirac delta 
distribution or point mass at the point $b$. All of this can be made 
rigorous, but at the expense of some long development in graduate level 
analysis. The approach here is at about the level of a second year course in 
engineering mathematics.

In addition, the Fourier transform behaves nicely on a product of 
functions and turns it into the convolution of the transforms:
\begin{equation}
   \widehat{fg} = \hat{f}\ast \hat{g}.
\end{equation}
For two functions $\phi(\omega)$ and $\psi(\omega)$, the convolution
$\phi \ast \psi$ is defined by
\begin{equation}
   (\phi \ast \psi)(\omega) = 
   \int_{-\infty}^{\infty}\phi(\alpha)\psi(\omega-\alpha) \, d\alpha.
\end{equation}
Again, we must extend convolution beyond the realm of functions. In 
particular we need convolutions of delta distributions and for them we 
can easily show that $\delta_0$ behaves as the identity for convolution
\begin{equation}
   \delta_0 \ast \phi = \phi
\end{equation}
and that
\begin{equation}
  \delta_a \ast \delta_b = \delta_{a+b}.
\end{equation} 

Now back to the cosine identity.  
Let $f(x)=\prod_{k=1}^{\infty} \cos (x/2^k)$ and let $f_n$ be the 
$n$th partial product. The Fourier transform of $f_n$ is
\[ \hat{f_n} = \ast \prod_{k=1}^n 
\frac{1}{2}(\delta_{1/{2^n}} + \delta_{-1/{2^n}}). \] 
The asterisk in front of the product sign indicates a repeated 
convolution of the factors. Expanding for $n=3$ we 
see that \[ \hat{f_3} = \frac{1}{8}(\delta_{-7/8} + \delta_{-5/8} + \cdots
\delta_{7/8}). \]
 Likewise \[\hat{f_n} = \frac{1}{2^n} \sum_{b \in B_n} 
\delta_b\, , \] where $B_n$ is the set of $2^n$ equally spaced numbers from
$-1 +  1/2^n$ to $1 - 1/2^n$ with spacing $ 2/2^n = 1/2^{(n-1)}$.

The sequence of measures $\hat{f_n}$ converges to the uniform density on $[-1,1]$
of total mass 1, which we can write as $(1/2)\chi_{[-1,1]} d\omega$.
The inverse transform is easy to compute: 
  \[ \int_{-1}^1 \frac{1}{2}e^{i\omega x} d\omega = \frac{\sin x}{x} .\]

The spectrum of $(\sin x)/x$ is uniform in the interval $ -1 \leq \omega 
\leq 1$. This means that $\sin x/x$ is a continuous linear combination
of the ``pure'' harmonics $e^{i\omega x}$ with the same weight of $1/2$
for each $\omega \in [-1,1]$. 

With this proof we have a way to generate a family of similar identities.
Let us put point masses at $3^n$ equally spaced points from $-1 + 
1/3^n$ to $1 - 1/3^n$ with spacing $2/3^n$. Such 
a measure is the convolution 
$\bigconv_{k=1}^n \frac{1}{3}(\delta_{-2/{3^k}} + \delta_0 
+\delta_{2/{3^k}})$.
Applying the inverse transform 
\[ {\cal F}^{-1}\left(\frac{1}{3}\left(\delta_{-2/{3^k}} + \delta_0 
   +\delta_{2/{3^k}}\right)\right) = 
   \frac{1}{3}\left(2\cos \frac{2x}{3^k} + 1\right) \]
 and taking 
limits gives us the infinite product identity
\begin{equation}
  \prod_{n=1}^\infty \frac{1}{3}\left(1 + 2 \cos \frac{2x}{3^n}\right) =
   \frac{\sin x}{x}.
\end{equation}
 
Let us use the positive integer $p$ as the base (we have just seen $p=2$ 
and $p=3$). The first measure $\hat{f_1}$ is the sum of point masses at 
$p$ points equally spaced from $-1+ 1/p $ to $1- 1/p$ with 
spacing $2/p$.
\begin{eqnarray}
 \hat{f_1} & = & \frac{1}{p}(\delta_{\frac{1-p}{p}} + \delta_{\frac{3-p}{p}}
  + \delta_{\frac{5-p}{p}} + \cdots + \delta_{\frac{p-1}{p}}) \\
           & = & \frac{1}{p}\sum_{l=0}^{p-1} \delta_{\frac{2l+1-p}{p}}
\end{eqnarray}
We let 
\[ \hat{f_n} = \bigconv_{k=1}^n \left( 
  \frac{1}{p}\sum_{l=0}^{p-1}\delta_{\frac{2l+1-p}{p^k}}\right)
\]
and one can see that $\hat{f_n}$ consists of $p^n$ point masses equally 
spaced from  $-1+1/p^n$ to $1-1/p^n$  with spacing 
$2/{p^n}$. Taking the inverse transform we see that
\[ {\cal F}^{-1}(\hat{f_n})(x) = \prod_{k=1}^n \frac{1}{p}\sum_{l=0}^{p-1}
  \exp((2l+1-p)ix/{p^k})
\]
Rewriting the exponentials as cosines and taking limits gives the general
identities. 

There is a slight difference in the form depending on the parity of $p$. 
For $p$ even 
\begin{equation} 
   \prod_{k=1}^\infty \frac{1}{p}\left(\sum_
   {\begin{array}{cc}
   1 \leq m \leq p-1 \\
   \mbox{$m$ odd}
   \end{array}}
   2\cos \frac{mx}{p^k} \right) 
   = \frac{\sin x}{x}. 
\end{equation}
For $p$ odd
\begin{equation}
   \prod_{k=1}^\infty \frac{1}{p}
   \left(1 + \sum_
   {\begin{array}{cc}
   1 \leq m \leq p-1 \\
   \mbox{$m$ even}
   \end{array}}2\cos \frac{mx}{p^k} 
   \right) = \frac{\sin x}{x}.
\end{equation} 
For $p=6$ the identity takes the form
\begin{equation}
   \prod_{k=1}^\infty \frac{1}{6}\left(2\cos\frac{x}{6^k}+2\cos\frac{3x}{6^k}+
    2\cos\frac{5x}{6^k}\right) = \frac{\sin x}{x}.  
\end{equation}
For $p=7$ the identity takes the form
\begin{equation}
   \prod_{k=1}^\infty \frac{1}{7}\left(1+2\cos\frac{2x}{7^k}+
    2\cos\frac{4x}{7^k} + 2\cos\frac{6x}{7^k}\right) = \frac{\sin x}{x}.  
\end{equation}
For larger $p$ fewer terms in the product are needed for the same degree 
of accuracy in the approximation to $\sin x/x$. In fact, by letting $p$ 
go to infinity the first factor alone approaches $\sin x/x$ and provides 
a novel derivation of a well-known result. I leave it to the reader to
work it out.
    
\section{Probabilistic Interpretation}

Mark Kac, in his delightful and now classic Carus monograph \cite{Kac59},
proves the first cosine identity (\ref{identity1}) in a way that 
is equivalent to the one 
we have outlined, although he does not explicitly
use the Fourier transform, delta functions, and convolution.
He then turns the identity into a question of probability, which for him 
was the leitmotif of his mathematical work.

The original product identity (\ref{identity1}) arises from the following 
experiment. Flip a fair coin repeatedly. Beginning with 0, add 
$1/2$ if the result is heads and subtract $-1/2$ if the result is 
tails. On the next toss add or subtract $1/4$; on the next add or 
subtract $1/8$, and so on. What is the distribution of the sums 
over the probability space whose elements are the countable sequences of 
coin tosses? Clearly, the sums are distributed uniformly between $-1$ and 
$1$.

Let $s_n$ denote the $n$th partial sum. It is a sum of independent random 
variables $a_1 + a_2 + \cdots + a_n$, where $a_k$ has the probability 
distribution $(1/2)\left( \delta_{1/{2^k}} + \delta_{-1/{2^k}} \right) $.
The probability distribution of a sum of independent random 
variables is the convolution of the respective distributions of the 
random variables. Therefore, $s_n$ has the distribution
\[ \bigconv_{k=1}^n \frac{1}{2}\left( \delta_{1/{2^k}} + \delta_{-1/{2^k}} 
  \right) = \frac{1}{2^n} \left( \delta_{\frac{1-2^n}{2^n}} + \cdots +
   \delta_{\frac{2^n-1}{2^n}} \right).  \]

The inverse Fourier transform of a probability measure is called its {\bf 
characteristic function}. Thus, the characteristic function for the 
distribution of $s_n$ is the product $\prod_{k=1}^n \cos x/2^k$. In 
the theory of probability and statistics, characteristic functions are a 
powerful tool. Typically computations are done with characteristic 
functions in order to draw conclusions about distributions of random 
variables as in the standard proof of the Central Limit Theorem. Here, 
however, we have inverted the relationship in order to compute with the 
probability measures and to get results about the characteristic functions.

\section{Related Products: Examples and Conjectures}

\subsection{Coin tossing and Cantor sets}

The Cantor set $K$ is the set of points between 0 and 1 whose ternary 
expansion has no 1's in it. So $z$ is in $K$ if   
$z=\sum_{k=1}^{\infty} t_k 3^{-k}$, $t_k \in \{0,2\}.$  Define $K_n$ to be 
the set of elements of $K$ that have the form $\sum_{k=1}^n t_k 3^{-k}$, and define
a probability measure supported on $K_n$ 
\begin{equation}
  \mu_n = \frac{1}{2^n}\sum_{z\in K_n}\delta_z .
\end{equation} 
$K_n$ has $2^n$ elements so $\mu_n$ is equally distributed on $K_n$.
The sequence $(\mu_n)$ has a limit $\mu$, which can be described as 
assigning the following limit as the measure of a set $E$:
\begin{equation} \mu(E) = \lim_{n\rightarrow \infty}\frac{\#E \cap K_n}{2^n}. 
  \label{Cantor.measure}
\end{equation}
The measure $\mu$ is also the Lebesgue-Stieltjes measure of the Cantor 
function. The Cantor function is continuous, non-decreasing, and has 
derivative zero on the complement of 
the Cantor set. Thus it defines a measure 
supported on the Cantor set, which is precisely the measure $\mu$ defined
in (\ref{Cantor.measure}).

What is of interest in this note is that 
$\mu_n$ is the finite convolution product
\be \mu_n= \bigconv_{k=1}^n \frac{1}{2}\left(\delta_0 + \delta_{2/{3^k}} \right).
\ee
Consider the experiment of tossing a fair coin. On toss number $k$ let
\[ a_k = \left\{ \begin{array}{cc} 0 & \mbox{heads} \\ 2/3^k & 
\mbox{tails} \end{array} \right.   \]
Let $s_n = \sum_{k=1}^{n}a_k$. Then $s_n$ is equally distributed over 
$K_n$. The characteristic function for the distribution of $s_n$ is
$\prod_{k=1}^{n}(1/2)\left(1 + e^{2xi/3^k} \right)$. Define 
\begin{equation} f(x) = \prod_{k=1}^{\infty}\frac{1}{2}\left(1 + e^{2xi/3^k} 
\right). \end{equation}
(One checks easily that the product is convergent.) Then $\hat{f}=\mu$, 
the Cantor measure, but is it possible to characterize $f$ in any other way? 

This leads us to look at the related infinite product 
$\prod_{k=1}^{\infty}\cos \frac{2x}{3^k}$. Because \[{\cal F}(\cos 
\frac{2x}{3^k}) = \frac{1}{2}\left(\delta_{2/3^k} + 
\delta_{-2/3^k}\right)\] the probabilistic interpretation is clear: add
or subtract $2/3^k$ on the $k$th toss with equal probability. Let $s_n$ 
be the sum of the first $n$ values. What is the distribution of $s_n$ and 
what is the distribution of $s=\lim_{n\rightarrow \infty} s_n$? The exercise of 
expanding and plotting the values of $s_3$ lead one to suspect that $s$ 
is distributed ``uniformly'' over the Cantor set constructed from 
$[-1,1]$ by successively removing middle thirds. That is easy to prove, 
as follows.

Define the affine map of $[0,1]$ to $[-1,1]$ by $z \mapsto 
2(z-1/2)$. Let $z=\sum t_k 3^{-k}$, $t_k \in \{0,2\}$, be a point in the 
Cantor set. The ternary expansion of $1/2$ is $\sum 3^{-k}$, and so
$2(z-1/2) = \sum 2(t_k - 1) 3^{-k}$. The coefficients $2(t_k - 1)$ are 
either 2 or $-2$ with equal probability.

This shows that the infinite product $\prod_{k=1}^{\infty}\cos \frac{2x}{3^k}$
has Fourier transform equal to the Cantor measure on the Cantor set
constructed from $[-1,1]$ by removing middle thirds, but it does not give 
us a closed form like $\sin x/x$. It would be most surprising if there were
any simpler description of $\prod_{k=1}^{\infty}\cos \frac{2x}{3^k}$. In Figure 2
is a plot of the partial product with $n=8$ and $0 \leq x \leq 100$. (The 
function is even.) Over this range the infinite product is 
indistinguishable from the eighth partial product. The self-similarity of 
the Cantor set at smaller and smaller scales appears to be reflected in 
the self-similarity of the graph at higher and higher frequencies.

\vspace{.4in}
\centerline {
\includegraphics[width=3in]{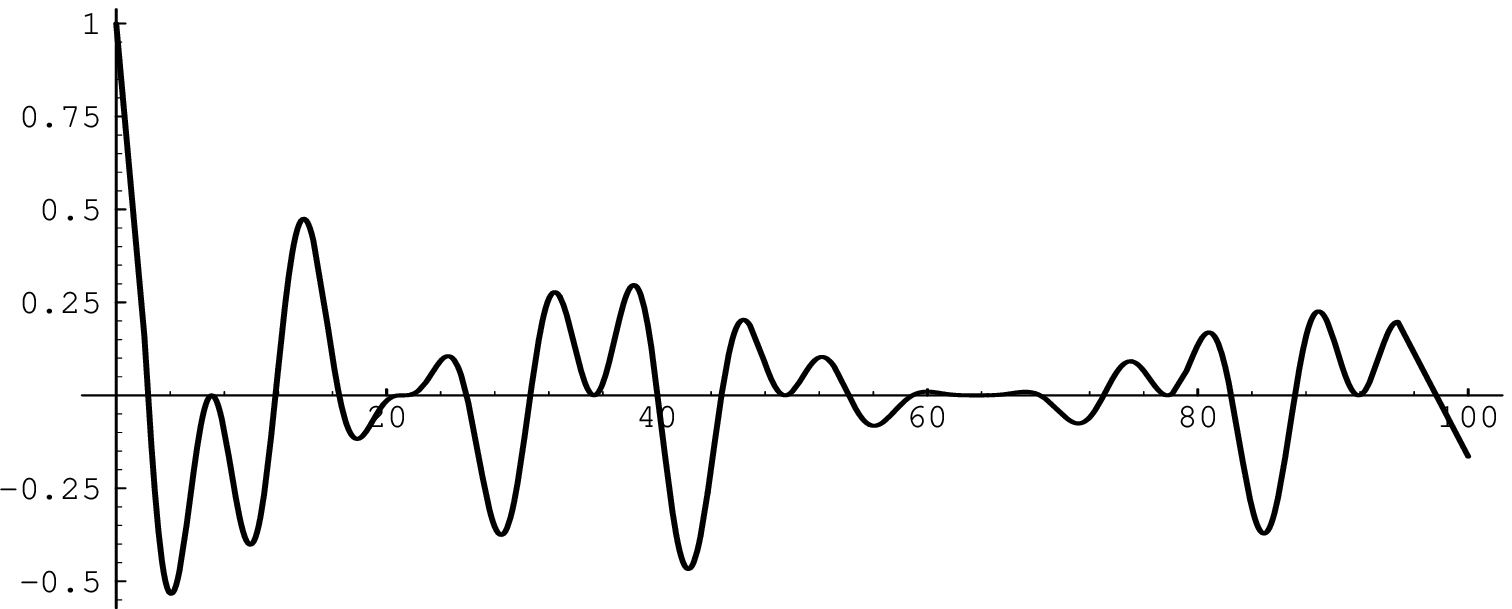}
}
\begin{center} 
 	Figure \thefigure.  
\end{center} 
\addtocounter{figure}{1}
    
\subsection{Harmonic Series With Random Signs}

We have been looking at the sums of series of the form
\be \sum_{k=1}^{\infty}t_k c_k  \ee
where $t_k$ is randomly chosen to be 1 or $-1$ with equal probability.
Rademacher proved that if $\sum {c_k}^2 < \infty$, then the sum converges 
with probability one on the probability space $\Omega = \{-1,1\}^{\N}$.
($\Omega$ can be identified with the unit interval and the probability 
measure with Lebesgue measure by using binary representations of numbers 
in the interval.) In \cite{Kac59} Kac gives the proof of this theorem due 
to Paley and 
Zygmund. It is also a theorem that the series diverges 
with probability one if $\sum {c_k}^2 
=\infty $. Let us consider the random harmonic series 
\be \sum_{k=1}^{\infty}\frac{t_k}{k}, \ee 
which converges almost surely by Rademacher's result,
with the goal of understanding the distribution of the sums.
This means we want to understand the distribution of the random variable 
$s$ defined on $\Omega$. If we let $s_n$ be the partial sum, also a 
random variable, then the probability distribution of $s_n$ is the measure
\be \mu_n = \bigconv_{k=1}^{n} \frac{1}{2}\left(\delta_{1/k} + 
\delta_{-1/k} \right)  \ee
and its inverse transform is
\be {\cal F}^{-1}(\mu)(x) = \prod_{k=1}^{n}\cos \frac{x}{k}. \ee
The product converges uniformly on compact sets as $n \rightarrow \infty$,
and so it is plausible that the sequence $\mu_n$ converges to a 
probability measure $\mu$ that is the distribution of the random variable $s$.
There is, however, a fair bit of analysis to make this rigorous. Assuming 
that the analysis can be made rigorous, then the plot of the Fourier 
transform of the infinite product $\prod_{k=1}^{\infty}\cos\frac{x}{k}$ 
will show how the sums are distributed. Let us call this function 
$\phi(\omega)$. Then 
\begin{eqnarray}
 \phi(\omega) & = & {\cal F}\left( \prod_{k=1}^{\infty}
    \cos\frac{x}{k} \right) (\omega) \\
    & = & \frac{1}{2\pi}\int_{-\infty}^{\infty}e^{-i\omega x}
    \prod_{k=1}^{\infty}\cos\frac{x}{k} \, dx \\
    & = & \frac{1}{2\pi}\int_{-\infty}^{\infty}(\cos \omega x + i\sin 
    \omega x) \prod_{k=1}^{\infty}\cos\frac{x}{k} \, dx \\ 
    & = & \frac{1}{\pi}\int_{0}^{\infty}\cos \omega x \prod_{k=1}^{\infty}
           \cos \frac{x}{k} \, dx .
\end{eqnarray}
There is not a closed form for $\phi(\omega)$ and so we resort 
to numerical integration. We truncated the infinite product at $n=1000$ and
integrated from 0 to 15 using a straightforward Riemann sum with 
$dx=0.02$ and the midpoints of the subintervals for the points of 
evaluation. Values for $\omega$ were from 0 to 3.8 in multiples of 0.2.
 The integration was done with True BASIC on a portable 
Macintosh. See Figure 3. The distribution is very flat for $-1 < \omega < 
1$, much flatter than a normal distribution. A few of the computed values 
are given in this table 
 
\renewcommand{\baselinestretch}{1} \large \normalsize
\begin{table} [h]
]\begin{center}
\begin{tabular}{r|r}
 \multicolumn{1}{c}{$\omega$} &
 \multicolumn{1}{c}{$\phi(\omega)$} \\ \hline
  0.0 & .249995 \\
  0.1 & .249991 \\
  0.2 & .249972 \\
  0.4 & .249809 \\
  0.6 & .249092 \\
  0.8 & .246819 \\
  1.0 & .241289 \\
  1.2 & .230494 \\
  1.4 & .212941 \\
  1.6 & .188425 \\
  1.8 & .158271 \\
  2.0 & .125000 \\
  2.2 & .091729 \\
  2.4 & .061576 \\
  2.6 & .030596 \\
  2.8 & .019506 \\
  3.0 & .008711 \\
  3.2 & .003181 \\
  3.4 & .000908 \\
  3.6 & .000192 \\
  3.8 & .000028
\end{tabular} \end{center}
\end{table}  
The value
of $\phi(0)$ is suspiciously close to $1/4$, suggesting perhaps that 
$\pi/4$ is the value of the integral
\begin{equation}
 \int_{0}^{\infty}\prod_{k=1}^{\infty}\cos \frac{x}{k} \, dx .         
\end{equation}
One might also conjecture that 
$\int_{0}^{\infty}\cos 2x\prod_{k=1}^{\infty}\cos \frac{x}{k} \, dx = 
\pi/8.$  \\
{\samepage\centerline {
\includegraphics[width=3in]{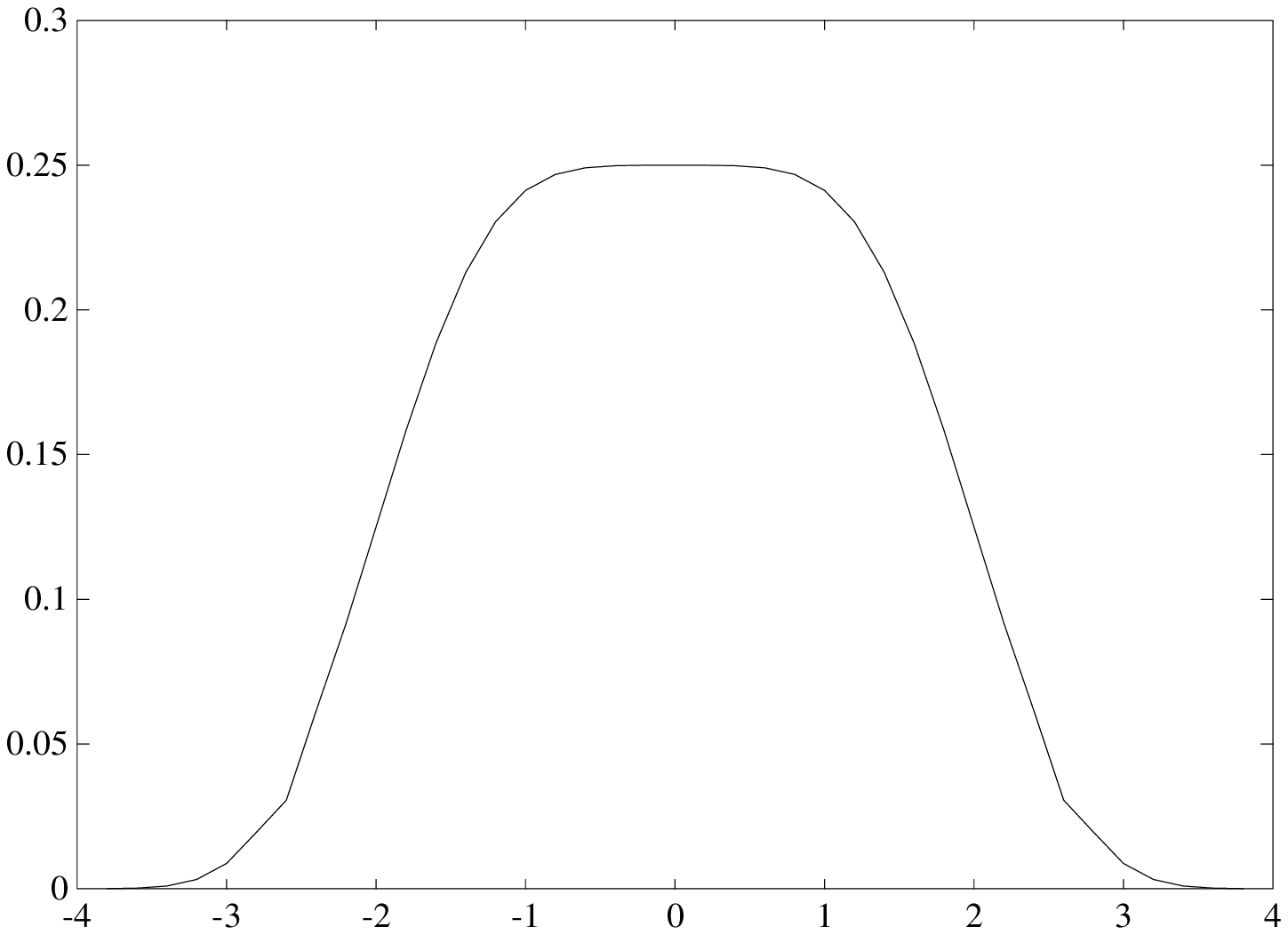}
}
  \begin{center} 
	  Figure \thefigure. Graph of $\phi(\omega)$.  
   \end{center}  }
   \addtocounter{figure}{1}
   
For additional evidence we turned to simulations of the sums.
Using MATLAB we ran 5000 sums of $\sum_{k=1}^{100} \frac{t_k}{k}$ with 
the values of $t_k$ picked randomly as $\pm 1$ with equal probability. 
Figure 4 shows a histogram of the sums.\\
\centerline {
\includegraphics[width=3in]{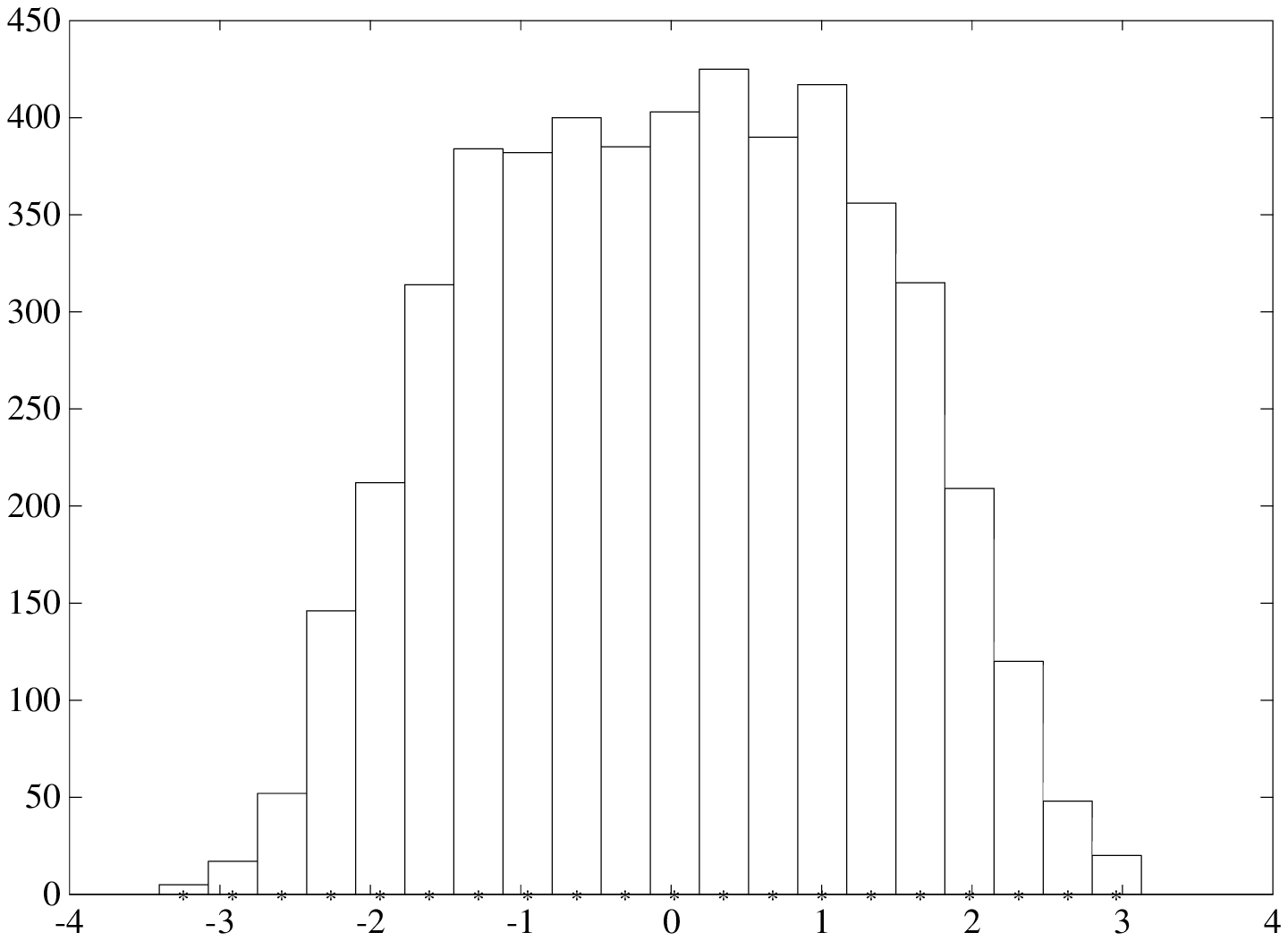}
}
	\begin{center} 
	  Figure \thefigure. Histogram of 5000 random sums.  
    \end{center}  \addtocounter{figure}{1}

\end{document}